\documentclass[12pt]{article}
\usepackage{theorem}
\usepackage{graphics}
\usepackage{amsfonts}
\usepackage{amssymb}

\theoremstyle{change}{\theorembodyfont{\slshape}
  \newtheorem{theorem}{Theorem.}[section]
  \newtheorem{proposition}[theorem]{Proposition.}
  \newtheorem{lemma}[theorem]{Lemma.}
  \newtheorem{remark}[theorem]{Remark.}
  \newtheorem{corollary}[theorem]{Corollary.}

    }

{ \theorembodyfont{\rmfamily} \newtheorem{example}[theorem]{Example.}
  
    }

\def\proof{\noindent{\bf Proof.}\enspace} \def\endproof{ \quad $\kasten$}

\def\tq{\mathop{\lower 4pt \hbox{$ {\buildrel{/} \over
        {\scriptstyle{\rm tor}}} $}}\nolimits} \def\cq{\mathop{\lower
    4pt \hbox{$ {\buildrel{/} \over {\scriptstyle{\rm cat}}}
      $}}\nolimits} \def\gq{\mathop{\lower 4pt \hbox{$ {\buildrel{/}
        \over {\scriptstyle{\rm good}}} $}}\nolimits}

\def\h#1{\widehat{#1}} \def\t#1{\widetilde{#1}}
\def\b#1{\overline{#1}}
\def\CC{{\mathbb C}} \def\ZZ{{\mathbb Z}} \def\RR{{\mathbb R}}
\def\NN{{\mathbb N}}  \def\PP{{\mathbb P}}

\def\S{\mathcal{S}}  \def\KMor{{\rm
    Mor}_{\mathfrak{K}}} \def\mal{\mathbin{\! \cdot \!}}

\def\cone{\mathop{\hbox{\rm cone}}}  \def\lin{\mathop{\hbox{\rm lin}}} \def\id{\mathop{\rm
    id}\nolimits} \def\pr{\mathop{\rm pr}\nolimits}
\def\Hom{\mathop{\rm Hom}\nolimits}  

\def\osubset{\subset \kern-10pt {\rm o}\ }

\def\topto#1{\mathop{\longrightarrow}\limits^{#1}}

\def\kasten{\mathord{\vbox{\hrule \hbox{\vrule \hskip5pt \vrule
        height5pt \vrule} \hrule}}}

\def\text#1{\hbox{\rm #1}}

\def\bigtopmapright#1{\smash{\mathop{\hbox to
      35pt{\rightarrowfill}}\limits^{#1}}}

\def\bigbotmapright#1{\smash{\mathop{\hbox to
      35pt{\rightarrowfill}}\limits_{#1}}}
\def\botmapright#1{\smash{\mathop{\hbox to
      30pt{\rightarrowfill}}\limits_{#1}}}

\def\bigtopmapleft#1{\smash{\mathop{\hbox to
      35pt{\leftarrowfill}}\limits^{#1}}}

\def\bigbotmapleft#1{\smash{\mathop{\hbox to
      35pt{\leftarrowfill}}\limits_{#1}}}

\def\rmapdown#1{\Big\downarrow\rlap{$\vcenter
    {\hbox{$\scriptstyle#1$}}$}} \def\lmapdown#1{\llap{$\vcenter{\hbox
      {$\scriptstyle#1$}}$}\Big\downarrow}

\def\rmapne#1{\nearrow\rlap{\hbox{$\scriptstyle#1$}}}

\def\rmapse#1{\searrow\rlap{\hbox{$\scriptstyle#1$}}}
\def\rmapsw#1{\swarrow\rlap{\hbox{$\scriptstyle#1$}}}

\def\lmapse#1{\llap{\hbox{$\scriptstyle#1$}}\searrow}

\begin{document}

\thispagestyle{empty}

\begin{center}

{\large\sc Examples and Counterexamples for}

\medskip

{\large\sc Existence of Categorical Quotients}

\bigskip

\bigskip

Annette A'Campo--Neuen %
and %
J\"urgen Hausen

%
%
%

\bigskip

\end{center}

{\small\noindent {\sc Abstract.}\enspace
We give examples for existence and non--existence of categorical
quotients for algebraic group actions in the categories of algebraic
varieties and prevarieties. All our examples are subtorus actions on
toric varieties.}

\section*{Introduction}

For group actions in the category of algebraic varieties,
various notions of quotients have been introduced. Among these, the
\emph{categorical quotient} is a basic concept; here
one only requires universality with respect to invariant morphisms. In
practice, it is a delicate problem whether or not a given action admits a
categorical quotient.
A possible way to obtain existence statements is to treat the
problem in a suitably  modified category. For example, if a finite
group acts on a variety, then this action in general admits no
algebraic variety as orbit space but in the category of algebraic
spaces it has a geometric quotient.

In this note we investigate the effect of allowing non--separated
quotient spaces on existence of categorical quotients. Our aim is to
show by means of examples that concerning categorical quotients the
separated and the non--separated case behave surprisingly independent from
each other.
We work with the following terminology:
Let $G$ be a complex algebraic group, let $\mathfrak{K}$ denote any
 subcategory of the category of complex algebraic prevarieties
containing $G$, and assume that $G$ acts $\mathfrak{K}$--morphically
on an object $X$ of $\mathfrak{K}$. Then we will call a morphism
$p\in\KMor(X,Y)$ a {\it $\mathfrak{K}$--quotient} for the action of $G$ on
$X$ if for every $G$--invariant morphism $f\in \KMor(X,Z)$ there is a
unique morphism $\t{f}\in\KMor(Y,Z)$ with $f=\t{f}\circ p$.

We consider actions of subtori on complex toric
varieties. In this setting, it makes
sense to ask for quotients in the categories $\rm{AV}$ of complex 
algebraic varieties, $\rm{PV}$ of complex algebraic prevarieties,
$\rm{TV}$ of complex toric varieties and $\rm{TP}$ of complex toric
prevarieties. In \cite{acha} and \cite{acha2} we have shown that
$\rm{TV}$-- and $\rm{TP}$--quotients always exist. For AV-- and
PV--quotients it is well--known that these notions need not coincide
if both exist (see Example \ref{hyperbolic}). Concerning existence and
non--existence we here give the following results and examples:

\begin{enumerate}
\item If $H$ is a subtorus of the big torus $T$ of a toric variety $X$
  with $\dim T/H \le 2$, then the $\rm{TV}$--quotient for the action
  of $H$ on $X$ is also an $\rm{AV}$--quotient.
\item A $\CC^*$--action without $\rm{AV}$--quotient but with
  $\rm{PV}$--quotient (see Section 5).
\item A $\CC^*$--action admitting neither an $\rm{AV}$--quotient nor a
  $\rm{PV}$--quotient (see Section~6).
\item A $\CC^*$--action with $\rm{AV}$--quotient and without
  $\rm{PV}$--quotient (see Section 7).
\end{enumerate}

The examples ii) and iii) in fact do not even admit a quotient in
the categories of algebraic or analytic spaces. The existence result
i) is proved in a slightly more general 
framework: Let $X$ be a toric prevariety and let $H$ be a subtorus of
the acting torus $T$ of $X$. We call a regular map $q \colon X \to Y$
an {\it $H$--invariant separation} of $X$, if $Y$ is a variety, $q$ is
$H$--invariant and every $H$--invariant regular map from $X$ to a
variety $Z$ factors uniquely through $q$. We prove:

\medskip

\noindent{\bf Theorem.}\enspace
{\sl If $\dim(T/{H}) \le 2$, then there exists an $H$--invariant
  separation of $X$.} 

\medskip

The present note is organized as follows: Sections 1 to 3 are devoted
to obtain some general criteria for existence and non--existence of
$H$--invariant separations and categorical quotients. The examples are
presented in Sections 4 to 7. Throughout the note we make use of the
basic concepts introduced in \cite{Fu} and \cite{acha}, \cite{acha2}.

\section*{Notation}

We fix some notation. Let $N$ be a lattice, i.e., a free $\ZZ$--module
of finite rank. The dual lattice is $M := \Hom(N,\ZZ)$. The
canonical pairing is denoted by
$$
M \times N \to \ZZ, \quad (u,v) \mapsto u(v) =: \langle u, v
\rangle .$$
Let $N_{\RR} := \RR \otimes_{\ZZ} N$ denote the real
vector space associated to $N$. Moreover, for a homomorphism $F \colon
N \to N'$ of lattices, denote by $F_{\RR}$ its extension to the real
vector spaces associated to $N$ and $N'$.

When we speak of a cone in $N$ we always think of a convex rational
polyhedral cone in $N_{\RR}$. For two cones $\tau$, $\sigma$ in $N$ we
write $\tau \prec \sigma$ if $\tau$ is a face of $\sigma$. The
relative interior of a cone $\sigma \subset N_{\RR}$ is denoted by
$\sigma^{\circ}$. The dual cone of a cone $\sigma$ in $N$ is the cone
$$\sigma^\vee := \{u \in M; \; \forall_{v \in \sigma} \ \langle u,v
\rangle \ge 0\}.$$
A {\it fan} in $N$ is a finite set $\Delta$ of strictly convex cones
in $N$ such that $\sigma, \sigma' \in \Delta$ implies $\sigma \cap
\sigma' \prec \sigma$ and $\sigma \in \Delta$ implies that also every
face of $\sigma$ lies in $\Delta$. For example for any given cone 
$\sigma$ the set
$\mathfrak{F}(\sigma)$ of its faces forms a fan in $N$. For two fans
$\Delta, \Delta'$ we will use the notation $\Delta\prec \Delta'$ if
$\Delta$ is a subfan of $\Delta'$.

A {\it system of fans} in $N$ is a finite family $\S :=
(\Delta_{ij})_{i,j \in I}$ of fans in $N$ such that $\Delta_{ij} =
\Delta_{ji}$ and $\Delta_{ij} \cap \Delta_{jk} \prec \Delta_{ik}$
holds for any $i,j,k \in I$. In particular, one has
 $\Delta_{ij}\prec \Delta_{ii}\cap\Delta_{jj}$ for all $i, j\in I$.
 A system
$(\Delta_{ij})_{i,j \in I}$ of fans is called {\it affine}, if for
every $i \in I$ the fan $\Delta_{ii}$ is the fan of faces of a single
cone $\sigma(i)$. The set of labelled cones of a system $\S =
(\Delta_{ij})_{i,j \in I}$ of fans is 
$$\mathfrak{F}(\S) :=
\{(\sigma,i); \; i \in I, \; \sigma \in \Delta_{ii}\}.$$
For a system $\S = (\Delta_{ij})_{i,j \in I}$ of fans in a lattice, we
define its {\it support} to be the set
$$ \vert \S \vert \; := \bigcup_{(\sigma,i) \in \mathfrak{F}(\S)}
\sigma . $$
In \cite{acha2} we showed that every affine system $\S$ of
fans defines a toric prevariety $X_{\S}$. Moreover, we introduced the
concept of a map of systems of fans and proved that the assignment $\S
\mapsto X_{\S}$ is an equivalence of categories.

\section{Factorization of Regular Maps}

In this section we prove a criterion for the existence of a
factorization of a regular map. The result may be of interest
independent from its application in our proof of
Theorem~\ref{invarsep}. By a {\it (pre--) variety} we mean throughout
this paper an algebraic (pre--) variety over the field $\CC$ of complex
numbers. Recall that any prevariety carries in a natural manner the
structure of a possibly non--hausdorff complex analytic space.

Let $X$ denote a prevariety. By a {\it local curve} in $x \in X$ we
mean a holomorphic mapping germ $\gamma \colon \CC_{0} \to X_{x}$ arising
from an algebraic curve, i.e., there is an algebraic curve $X'$ in
$X$ through $x$ with $\gamma(\CC_{0})\subset X'_{x}$. Let
$p \colon X \to Y$ be a regular map of prevarieties. We say that a
local curve $\t{\gamma}$ in $x \in X$ is a {\it weak $p$--lifting} of
a local curve $\gamma$ in $y \in Y$ if there is a non--constant
holomorphic mapping germ $\alpha \colon \CC_{0} \to \CC_{0}$ and a
commutative diagram
$$
\matrix{ \CC_{0} & \bigtopmapright{\t{\gamma}} & X_{x} & \cr
  \lmapdown{\alpha} & & \rmapdown{p} & \cr \CC_{0} &
  \bigtopmapright{\gamma} & Y_{y} & .}$$
We call the map $p$ {\it
  weakly proper}, if any local curve in $Y$ admits a weak
$p$--lifting. Note that a weakly proper map is necessarily
surjective. Moreover, every proper regular map is weakly proper. For
a related notion in the context of algebraic spaces see \cite{Ko},
Section~3.

\begin{proposition}\label{factmor}
  Let $p \colon X \to Y$ be a weakly proper regular map of
  prevarieties, assume that $Y$ is normal and let $f \colon X \to Z$
  be a regular map into a variety. If $f$ is constant on the fibres of
  $p$, then there is a unique regular map $\t{f} \colon Y \to Z$ such
  that $f = \t{f} \circ p$.
\end{proposition}

\proof By assumption $\t{f}$ exists as a uniquely determined map of
sets. We have to show that $\t{f}$ is regular. Since $Y$ was assumed
to be normal, it suffices to show that the graph $\Gamma$ of $\t{f}$ is
closed with respect to the Zariski topology in $Y\times Z$. Note that
$$\Gamma = \{(p(x),f(x)); \; x \in X\}.$$
Hence $\Gamma$ is a constructible subset of $Y \times Z$. In
particular the Zariski closure $\b{\Gamma}$ of $\Gamma$ in $Y\times
Z$ coincides with the metric closure of $\Gamma$ and there is a dense
subset $\Gamma^{0} \subset \Gamma$ which is Zariski open in
$\b{\Gamma}$. Now, let $(y,z)$ be a point of $\overline{\Gamma}
\subset Y \times Z$.

Since $\b{\Gamma} \setminus \Gamma^{0}$ is nowhere dense and
Zariski closed, we can choose an open disc $U \subset \CC$ around $0$
and a holomorphic curve $\gamma \colon U \to \overline{\Gamma}$ such
that the Zariski--closure of $\gamma(U)$ is an algebraic curve,
$\gamma(0) = (y,z)$ holds and the complement of $U_1 :=
\gamma^{-1}(\Gamma)$ in $U$ is discrete and closed (see
e.g. \cite{Kr}). Since $p$ is weakly proper, we find an open
neighbourhood $V$ of $0 \in \CC$, a regular curve $\t{\gamma} \colon V
\to X$ and a non--constant regular map $\alpha \colon V \to U$ such that
$$\alpha(0) = 0, \qquad p \circ \t{\gamma} = \pr_{Y} \circ \gamma
\circ \alpha. $$
Since $\alpha$ is non--constant, $\alpha^{-1}(U_1)$
has a closed discrete complement in $V$. For any $s \in
\alpha^{-1}(U_1)$ we have
$$
\gamma(\alpha(s)) = (p(\t{\gamma}(s)), f(\t{\gamma}(s))). $$
Thus, for continuity reasons, we obtain
$$
(y,z) = \gamma(\alpha(0)) = (p(\t{\gamma}(0)), f(\t{\gamma}(0)))
\in \Gamma. \qquad \kasten$$

\medskip

For open surjections $p$ the above result is well--known (see
e.g. \cite{Bo}, II.6.2). We will apply Proposition~\ref{factmor} in the
following situation. Let $\S$ be a system of fans in a lattice $N'$,
let $\Delta$ be a fan in a lattice $N$ and assume that $P \colon N'
\to N$ is surjective and defines a map of systems of fans from $\S$ to
$\Delta$. Denote by $p \colon X_{\S} \to X_{\Delta}$ the toric
morphism associated to $P$. Then we obtain the following
characterization of weak properness:

\begin{proposition}\label{weakpropchar}
  The map $p$ is weakly proper if and only if $P_{\RR}(\vert \S \vert)
  = \vert \Delta \vert$.
\end{proposition}

For the proof we formulate some auxiliary results. Let $T$ be the
acting torus of $X_{\Delta}$ and denote by $x_{0}$ the base point of
$X_{\Delta}$. Call a local curve in $X_{\Delta}$ {\it generic} if its
image intersects the open orbit $T \mal x_{0}$. For any $v \in N$, we
denote by $\lambda_{v} \colon \CC^{*} \to T$ the associated
one--parameter--subgroup.

\begin{lemma}\label{curvedecomp}
  Let $\gamma$ be a generic local curve in $X_{\Delta}$. Then there
  exists a local curve $\beta$ in $t \in T$ and a point $v \in \vert
  \Delta \vert \cap N$ such that near $0$ we have $\gamma(s) =
  \beta(s)\lambda_{v}(s) \mal x_{0}.$
\end{lemma}

\proof We may assume that $N = \ZZ^{n}$ and hence $T = (\CC^{*})^{n}$.
Let $\gamma$ be defined on some open disc $U \subset \CC$ around
$0$. Set $V := \gamma^{-1}(T \mal x_{0})$. Then $U \setminus V$ is a
proper analytic subset of $U$ and hence it is discrete and
closed. Thus, after shrinking $U$, we may assume that either $U = V$
or $U = V \cup \{0\}$ holds. On $V$ there is a representation
$$
\gamma(s) = (g_{1}(s), \ldots, g_{n}(s)) \mal x_{0}$$
with holomorphic functions $g_{i} \in \mathcal{O}_{\rm
  an}^{*}(V)$. Using Laurent series expansion, we obtain $g_{i}(s) =
s^{v_{i}} \beta_{i}(s)$ with an integer $v_{i}$ and a function
$\beta_{i} \in \mathcal{O}_{\rm an}^{*}(U)$. Let $v := (v_{1},
\ldots, v_{n})$ and $\beta := (\beta_{1}, \ldots, \beta_{n})$. Then
$$
\lim_{s \to 0} \lambda_{v}(s) \mal x_{0} = \beta(0)^{-1} \mal
\gamma(0). $$
Consequently, the point $v$ lies in $\vert \Delta \vert$
and the desired decomposition of $\gamma$ is given by $\gamma =
\beta(s) \lambda_{v}(s) \mal x_{0}$. \endproof

\medskip

Now, let $\Delta_{1}, \ldots, \Delta_{r}$ be fans in lattices $N'_{i}$
and let surjective $P_{i} \colon N_{i}' \to N$ be given that are maps
of the fans $\Delta_{i}$ and $\Delta$. Let $p_{i} \colon
X_{\Delta_{i}} \to X_{\Delta}$ be the associated toric morphisms.

\begin{lemma}\label{famoffans}
  Let $\gamma$ be a generic local curve in $X_{\Delta}$. If $ \vert
  \Delta \vert = P_{\RR}(\vert \Delta_{1} \vert) \cup \ldots \cup
  P_{\RR}(\vert \Delta_{r} \vert)$, then, for some $i$, there is a
  weak $p_{i}$--lifting of $\gamma$.
\end{lemma}

\proof Choose $\beta$ and $v \in N \cap \vert \Delta \vert$ as in
Lemma~\ref{curvedecomp}. By assumption, for some $i$, there is a $v'
\in \vert \Delta_{i} \vert$ such that $P_{i}(v') = lv$ with a positive
integer $l$. Moreover, since $P_{i}$ is surjective, we have a
splitting
$$
\matrix{\quad T'_{i} \quad \topto{\cong}\ T \times \ker(\pi_{i}) &
  \cr \lmapse{\pi_{i}} \qquad \rmapsw{\pr_{T}} & \cr T &,}$$
where
$T'_{i}$ is the acting torus of $X_{\Delta_{i}}$ and $\pi_{i} \colon
T_{i}' \to T$ is the homomorphism associated to $p_{i}$. In
particular, there is a lifting $\t{\beta}$ with respect to $\pi_{i}$
of the local curve $s \mapsto \beta(s^{l})$ in $t \in T$. Now, let
$x_{0}'$ be the base point of $X_{\Delta_{i}}$. Then the desired weak
$p_{i}$--lifting of $\gamma$ is given by
$$
\t{\gamma} (s) := \t{\beta}(s) \lambda_{v'}(s) \mal x_{0}', \qquad
\alpha(s) := s^{l}. \qquad \kasten$$

\medskip

Finally, we need an elementary statement from convex geometry. Let
$\sigma$ denote a strictly convex polyhedral cone in some real vector
space, let $\tau$ be a face of $\sigma$ and let $P \colon V \to
V/\lin(\tau)$ be the projection.

\begin{lemma}\label{conecover}
  If $\sigma = \sigma_{1} \cup \ldots \cup \sigma_{r}$ with polyhedral
  cones $\sigma_{i}$, then $P(\sigma)$ is the union of all
  $P(\sigma_{i})$, where $\tau^{\circ} \cap \sigma_{i} \ne \emptyset$.
\end{lemma}

\proof We prove the assertion by induction on $r$. For $r=1$ or $\tau
=\{0\}$ there is nothing to show, so assume that $r>1$ and $\tau$ is
not trivial. Suppose that for some $j$ and some $v\in\sigma_j$ we had
$P(v) \notin P(\sigma_i)$ for all $\sigma_i$ meeting $\tau^{\circ}$.
Then $v$ does not lie in $\tau$. Hence there is a linear form
$u\in\sigma_j^{\vee}$ with $u(v)>0$ and $u(w)<0$ for some
$w\in\tau^{\circ}$. Fix a large $n\in\NN$ such that $u(v+nw)<0$.  Now
consider the cone
$$\sigma':=\sigma \cap \{v\in N_{\RR}; \; u(v)\leq 0\}.$$
Note that
$\sigma'$ contains $v+nw$ and is covered by less then $r$ of the cones
$\sigma_{i}' := \sigma_{i} \cap \sigma$. Moreover,
$\tau':=\sigma'\cap\tau$ is a face of $\sigma'$ and has the same
dimension as $\tau$. The induction hypothesis provides an $i$
and a $w'\in\lin(\tau') = \lin(\tau)$ such that $v+nw+w'\in
\sigma_i'$ and $\emptyset\ne\sigma_i'\cap
(\tau')^{\circ}\subset\tau^{\circ}$. \endproof

\medskip

\noindent
{\bf Proof of Proposition \ref{weakpropchar}.}\enspace First assume
that $p$ is weakly proper. Clearly $P_{\RR}(\vert \S \vert) \subset
\vert \Delta \vert$. To obtain the reverse inclusion, assume that
there are points $v \in \vert \Delta \vert \setminus P_{\RR}(\vert \S
\vert)$. Then we even find such a $v$ lying in $N$. For this $v$, the
curve $\lambda_{v}(s) \mal x_{0}$ admits no weak $p$--lifting,
contradicting our assumption on $p$.

Now, assume that $P_{\RR}(\vert \S \vert)$ equals $\vert \Delta
\vert$. Let $U \subset \CC$ be an open disc around zero and let $\gamma
\colon U \to X_{\Delta}$ be a holomorphic curve such that the 
 Zariski--closure of $\gamma(U)$ is an algebraic curve. Then there is
 a unique $T$--orbit $T \mal x_{\tau}$ of minimal dimension such that
 $\gamma(U)$ is contained in $V_{\tau} := \b{T \mal x_{\tau}}$.  Note
 that $V_{\tau}$ is itself a toric variety.

We will use the fact that $p^{-1}(V_{\tau})$ is a union of toric
varieties to apply Lemma \ref{famoffans}. Let $x_{\sigma} \in
V_{\tau}$, i.e., $\sigma$ is a cone of $\Delta$ with $\tau \prec
\sigma$. By our assumption, we have $\sigma = P_{\RR}(\sigma_{1}) \cup
\ldots \cup P_{\RR}(\sigma_{s})$ with certain $(\sigma_{i},k_{i}) \in
\mathfrak{F}(\S)$. By suitable ordering we achieve that
$P_{\RR}(\sigma_{i})$ meets $\tau^{\circ}$ if and only if $i \le r$
with some $r \le s$. Now, the above Lemma \ref{conecover} implies
$$\sigma + \lin(\tau) \; = \; \bigcup_{i=1}^{r} P_{\RR}(\sigma_{i}) +
\lin(\tau) \; \subset \; N_{\RR} / \lin(\tau). $$
Set $\tau_{i} :=
P_{\RR}^{-1}(\tau) \cap \sigma_{i}$ and consider the orbit closures
$V_{\tau_{i}}$ in $X_{\sigma_{i}}$. Note that $p(x_{\tau_i})=x_{\tau}$
since $P_{\RR}(\tau_i)\cap\tau^{\circ}\ne \emptyset$.  Therefore $p$
induces toric morphisms $p_{i} \colon V_{\tau_{i}} \to V_{\tau}$.
Applying this procedure to all the other $\sigma \in \Delta$ with
$\tau \prec \sigma$, we obtain a family of locally closed toric
varieties $V_{\tau_{j}} \subset p^{-1}(V_{\tau})$ and toric morphisms
$V_{\tau_{j}} \to V_{\tau}$. According to \cite{acha}, Example~2.7
these toric morphisms satisfy the assumptions of
Lemma~\ref{famoffans}. \endproof

\section{Two Cones}

In this section we consider the special case of a toric prevariety $X$
arising from an affine system $\mathcal{S}$ of fans in a lattice $N$
with two maximal cones $\sigma(1)$ and $\sigma(2)$. Let $L$ be a
primitive sublattice of $N$. Throughout this section we assume that
the projection $P \colon N \to N/L$ satisfies
$$
P_{\RR}(\sigma(1))^{\circ} \cap P_{\RR}(\sigma(2))^{\circ} \ne
\emptyset. \leqno{(*)}$$
Let $H$ be the subtorus of the big torus $T$
of $X$ corresponding to $L \subset N$ and suppose that $f\colon X\to
Z$ is an $H$--invariant regular map to a (not necessarily toric)
variety $Z$. A first simple observation is

\begin{lemma}\label{xsigmaigleichewerte}
  Let $t\in T$. Then we have:
\begin{enumerate}
\item There are regular curves $C_1, C_2 \colon \CC \to X$ and $C
  \colon \CC^* \to H$ with $C_1(s) = C(s) C_2(s)$ for all $s \in
  \CC^*$ and $C_i(0) = t \mal x_{[\sigma(i),i]}$.
\item $f(t \mal x_{[\sigma(1),1]}) = f(t \mal x_{[\sigma(2),2]})$. In
  particular, $f$ is constant on $T' \mal x_{[\sigma(1),1]}$, where
  $T' := T_{x_{[\sigma(1),1]}} T_{x_{[\sigma(2),2]}}$.
\end{enumerate}
\end{lemma}

\proof By assumption $(*)$, there are $w_i \in \sigma(i)^{\circ} \cap
N$ such that $w_{1} = v_{L} + w_{2}$ holds for some $v_{L} \in L$. Let
$\lambda_{w_i}$ and $\lambda_{v_L}$ denote the
one--parameter--subgroups of $T$ corresponding to these lattice
vectors. The curves
$$C_i \colon \CC^* \to X, \quad s \mapsto t \mal \lambda_{w_i}(s) \mal
x_0 $$
can be extended regularly to $\CC$ by setting $C_i(0) := t \mal
x_{[\sigma(i),i]}$. Together with the curve $C \colon \CC^* \to H$, $s
\mapsto \lambda_{v_L}(s)$ the $C_i$ satisfy i).

In order to check ii) note that according to i) the points
$x_{[\sigma(1),1]}$ and $x_{[\sigma(2),2]}$ cannot be separated by
$H$--stable complex open neighbourhoods. Since $f$ is continuous with
respect to the complex topology and $Z$ is hausdorff, the claim
follows. \endproof

\medskip

Note that in the proof of assertion ii), we only used that $f$ is
$H$--invariant and continuous with respect to the complex
topology. Hence the statement holds also for holomorphic $f$.

\begin{proposition}\label{zweikegel}
  Assume that there are faces $\tau_{i} \prec \sigma(i)$, and $v_i \in
  \tau_i^{\circ} \cap N$ such that the cone generated by $P(v_1)$ and
  $P(v_2)$ is a line. Then $f(\lambda_{v_i}(\CC^*)\mal x)=f(x)$ for
  $i=1,2$ and for all $x\in X$.
\end{proposition}

\proof Since $X$ is covered by affine $T$--stable open subspaces, it
suffices to show that with some nonempty open subset $V$ of $T$ we
have $f(\lambda_{v_i}(s) \mal t \mal x_0)=f(t \mal
x_0)$ for all $s\in\CC^*$ and all $t$ contained in $V$.

By appropriate scaling we achieve that $v_1 + v_2 \in L$ holds.
Assumption $(*)$ provides $w_i \in \sigma(i)^{\circ} \cap N$ such that
$w_1-w_2\in L$. Let $X_i$ denote the affine chart of $X$ corresponding
to $\sigma(i)$.  We now want to define
toric morphisms $\varphi_i\colon \CC\times\CC\times T\to X_i$.

We consider the
lattice homomorphisms $F_i\colon \ZZ^2\times N \to N$,
defined by $F_i(e_1)=v_i$, $F_i(e_2)=w_i$, and
$F_i(v)=v$ for all $v\in N$.
The corresponding toric morphisms are the following maps:
$$\varphi_i \colon \CC \times \CC \times T \to X_i, \quad
(s,r,t)\mapsto \cases{ t \mal \lambda_{v_i}(s) \mal \lambda_{w_i}(r)
  \mal x_0 & if $r \ne 0 \ne s$, \cr t \mal \lambda_{w_i}(r) \mal
  x_{[\tau_i,i]} & if $r\ne 0$, $s=0$, \cr t \mal x_{[\sigma(i),i]} &
  if $r=0$. \cr }$$

We use the regular maps $\varphi_i$ to define a regular map $\psi
\colon \CC \times \PP_1 \times T \to Z$. First notice that
$H$--invariance of $f$ yields for all $s \in \CC^*$, $r \in \CC$ and
$t \in T$ the identity
$$f(\varphi_1(s,r,t))= f(t\mal \lambda_{v_1}(s)\mal
\lambda_{w_1}(r)\mal x_0)= f(t\mal \lambda_{v_2}(1/s)\mal
\lambda_{w_2}(r)\mal x_0)= f(\varphi_2(1/s,r,t))\,.$$
So the rational
map $\PP_1 \times \CC \times T \to Z$ given by $([s_0,s_1],r,t)
\mapsto f(\varphi_1(s_1/s_0,r,t))$ extends to a morphism $\psi$. The
fibre $\psi^{-1}(z)$ of $z := \psi(0,0,e_{T})$ contains $\PP_1 \times
\{(0,e_T)\}$, where $e_T$ denotes the neutral element of $T$.

Now choose an open affine neighbourhood $W$ of $z$ in $Z$ and set $Y
:= \PP_1 \times \CC \times T \setminus \psi^{-1}(W)$. Consider the
projection $\pr \colon \PP_1 \times \CC \times T \to \CC \times T$.
Since $\PP_1$ is complete, $\pr(Y)$ is closed in $\CC \times T$.
Moreover, we have $\pr(z) \not\in \pr(Y)$. Thus for $W_0 := \CC \times
T \setminus \pr(Y)$ we have
$$
\psi^{-1}(z) \subset \PP_1 \times W_0 \; = \;\PP_1 \times \CC
\times T \setminus \pr^{-1}(\pr(Y)) \; \subset \; \psi^{-1}(W). $$
Since we chose $W$ to be affine, $\psi$ maps $\PP_1 \times \{w\}$ to a
point for every $w \in W_{0}$. In particular, for every point $(r,t)
\in W_0 \cap \CC^* \times T$ we have
$$
f(\lambda_{w_1}(r) \mal \lambda_{v_1}(s) \mal t\mal x_0)=
f(\lambda_{w_1}(r) \mal t \mal x_0)$$
for all $s \in \CC^*$. So $f$ is
constant on orbits of the one parameter subgroup $\lambda_{v_1}$ on a
dense subset of $T \mal x_0 \subset X_{i}$ and hence this is true
everywhere.  Since $v_1 + v_2 \in L$, this also holds for
$\lambda_{v_2}$.  \endproof

\section{A Criterion for the Existence of an Invariant Separation}

Let $\mathcal{S} = (\Delta_{ij})_{i,j\in I}$ denote a system of fans
in a lattice $N$ and let $L$ be a sublattice of $N$. Suppose that the
projection $P \colon N \to N/L =: \t{N}$ fullfills the following
conditions:
\begin{enumerate}
\item $\tau := \bigcup_{i \in I} P_{\RR}(\sigma(i))$
  is a strictly convex cone.
\item For every face $\varrho \prec \tau$ and any two $(\sigma,i),
  (\sigma',i') \in \mathfrak{F}(\S)$ with
  $P_{\RR}(\sigma)^{\circ} \cup P_{\RR}(\sigma')^{\circ} \subset
  \varrho^{\circ}$, there is a chain
$$ (\sigma,i)=(\sigma_{i_1},i_1), \ldots, (\sigma_{i_r},i_r) =
  (\sigma',i')$$
 in $\mathfrak{F}(\S)$ such that each
 $P_{\RR}(\sigma_{i_{k}})^{\circ}$ is contained in $\varrho$ and
 $P_{\RR}(\sigma_{i_{k}})^{\circ} \cap
 P_{\RR}(\sigma_{i_{k+1}})^{\circ} \ne \emptyset.$
\end{enumerate}

\begin{remark}
  If $\dim(\t{N}) \le 1$ then i) implies ii).  If $\dim\t{N}
  = 2$, then ii) is equivalent to $\tau^{\circ} =
  \bigcup_{\dim(P_{\RR}(\sigma(i))=2} P_{\RR}(\sigma(i))^{\circ}$. \endproof
\end{remark}

The projection $P$ defines a map of systems of fans from $\S$ to the
fan of faces of $\tau$. Moreover, denoting by $H$ the subtorus of $T$
that corresponds to $L$, we have

\begin{proposition}\label{sepcrit}
  The toric morphism $p \colon X_{\S} \to X_{\tau}$ defined by $P$ is
  an $H$-invariant separation.
\end{proposition}

\proof By Propositions \ref{factmor} and \ref{weakpropchar} it
suffices to show that every $H$--invariant morphism $f \colon
X_{\mathcal{S}} \to Z$ to a variety is constant on the fibres of
$p$. So let $\pi \colon T \to \t{T}$ denote the homomorphism of the
acting tori associated to $p$. Then the $p$--fibre of a point $\t{t}
\mal x_{\varrho} \in X_{\tau}$ is
$$
p^{-1}(\t{t}\cdot x_{\varrho}) = \bigcup_{P_{\RR}(\sigma)^{\circ}
  \subset \varrho^{\circ}} \pi^{-1}(\t{t} \mal \t{T}_{x_{\varrho}})
\mal x_{[\sigma,i]}$$
(see \cite{acha2}, Proposition 3.5). Let $T'$
denote the subtorus of $T$, generated by all isotropy groups
$T_{x_{[\sigma,i]}}$, where $P_{\RR}(\sigma^{\circ}) \subset
\varrho^{\circ}$, i.e., $T'$ correpsonds to the maximal sublattice in
the vector subspace spanned by the $\lin(\sigma)$. Then $T'\cdot
H=\pi^{-1}(\t{T}_{x_{\varrho}})$.

Now, for $(\sigma, i)$ and $(\sigma',i') \in \mathfrak{F}(\S)$ with
$P_{\RR}(\sigma^{\circ})\cup
P_{\RR}({\sigma'}^{\circ})\subset\varrho^{\circ}$, the chain condition
ii) implies by Lemma \ref{xsigmaigleichewerte} that $ f(t \mal
x_{[\sigma,i]}) = f(t \mal x_{[\sigma',i']})$ for all $t$, and hence
that $f$ is constant on $T' \mal t \mal x_{[\sigma,i]}$.  That shows
that $f$ is constant on the fibre $p^{-1}(\t{t} \mal x_{\varrho})$.
\endproof

\section{Codimension Two}

Let $X$ be a toric prevariety and let $H$ be a subtorus of the acting
torus $T$ of $X$. Denote by $\h{H}$ the maximal subtorus of $T$ such
that every $H$--invariant regular map from $X$ to a variety $Z$ is
invariant by $\h{H}$. In this section we prove

\begin{theorem}\label{invarsep}
  If $\h{H}$ is of codimension at most two in $T$, there exists an
  $H$--invariant separation for $X$.
\end{theorem} 

\begin{corollary}\label{presurfsep}
  Every toric presurface admits a separation. \endproof
\end{corollary}

\begin{corollary}\label{codimtwo}
  If $X$ is a toric variety and $H$ is of codimension at most two in
  $T$, then the TV--quotient is also an AV--quotient.  \endproof
\end{corollary}

As we shall see in Section 7, a PV--quotient need not exist even in small
codimension, and even if it exists, the AV--quotient and the
PV--quotient may be different (see Example 4.6).

For the proof of Theorem \ref{invarsep}, we may assume that $H =
\h{H}$. In particular, $H$ itself is of codimension at most two in
$T$. Moreover, we may assume that $X = X_{\S}$ for some affine system
of fans $\S$ in a lattice $N$.

Let $L$ denote the sublattice of $N$ that corresponds to $H$ and let
$P \colon N \to N/L =: \t{N}$ be the projection. Define an equivalence
relation on the index set $I$ by
$$
i \sim j :\iff \exists \ i= i_1, \ldots, i_r = j \hbox{ with }
P_{\RR}(\sigma(i_k))^{\circ} \cap P_{\RR}(\sigma(i_{k+1}))^{\circ} \ne
\emptyset .$$

\begin{lemma}
  For each equivalence class $E \subset I$ the set $\tau_E :=
  \bigcup_{\sigma \in E} P_{\RR}(\sigma(i))$ is a strictly convex
  cone.
\end{lemma}

\proof Maximality of $H$ and Proposition \ref{zweikegel} imply that
every cone $P_{\RR}(\sigma(i))$ is strictly convex. In particular the
assertion is verified in the case $\dim(T/H) \le 1$. Now suppose that
$\t{N}$ is of dimension two and there is an equivalence class $E$ such
that $\tau_E$ is not strictly convex. Then we find subsets $E_1$,
$E_2$ of $E$ such that each
$$
\tau_k := \bigcup_{\sigma \in E_k} P_{\RR}(\sigma)$$
is strictly
convex, $\tau_1^\circ \cap \tau_2^\circ \ne \emptyset$ and $\tau_1
\cup \tau_2$ is not strictly convex. Let $X_k$ be the open $T$-stable
subspace of $X_{\mathcal{S}}$ defined by the cones $\sigma(i)$ with
$P_{\RR}(\sigma(i)) \subset \tau_k$ According to Proposition
\ref{sepcrit}, the map $P$ defines a $H$-invariant separations $p_k
\colon X_k \to X_{\tau_k}$, $k=1,2$.

Now let $f \colon X \to Z$ be any $H$-invariant regular map. Set $f_k
:= f \vert_{X_k}$. Then we obtain the following commutative diagram of
regular maps:

$$ \matrix{%
X_1 \cup_T X_2 \quad \bigtopmapright{f_1 \cup_T f_2} \quad Z &\cr %
\lmapse{p_1 \cup_T p_2} \qquad \qquad \rmapne{\t{f}} &\cr%
X_{\tau_1} \cup_{\t{T}} X_{\tau_2} & .%
}$$

Here ``$\cup_T$'' indicates glueing along $T$. Let $\t{L} \subset
\t{N}$ be a line contained in $\tau_1 \cup \tau_2$. Then Proposition
\ref{zweikegel} yields that $\t{f}$ is invariant with respect to the
action of the subtorus $\t{H} \subset \t{T}$ corresponding to $\t{L}$.
Now let $\pi \colon T \to \t{T}$ denote the homomorphism of the acting
tori determined by $p$. Then $f_1 \cup_T f_2$ and hence $f$ is
invariant by $\pi^{-1}(\t{H})$. This contradicts the maximality of
$H$. \endproof

\medskip

\noindent
{\bf Proof of Theorem \ref{invarsep}.}\enspace By construction, the
cones $\tau_E$, where $E$ runs through the equivalence classes of
$\sim$, form a fan $\Delta$ in $\t{N}$. Moreover, the projection $P$
determines a map of systems of fans from $\S$ to $\Delta$. It follows
directly from Proposition~\ref{sepcrit} that the
associated toric morphism $p$ is a $H$--invariant separation of
$X$. \endproof

\begin{example}\label{hyperbolic}
  {\it A $\CC^{*}$--action with AV-- and PV--quotient different from
    each other.}  Let $X := \CC^{2} \setminus \{ 0\}$ and consider the
  action of $H := \{(t,t^{-1}); t \in \CC^{*}\} \subset (\CC^{*})^{2}$
  on $X$. Then the AV--quotient for this action is by
  \ref{codimtwo} just the map
  $$
  X \to \CC, \qquad (z,w) \mapsto zw. $$
  On the other hand, the
  PV--quotient is given by the following map from $X$ onto
  the line $\CC_{00}$ with doubled zero
  $$
  (z,w) \mapsto \cases{zw & if $zw \ne 0$, \cr 0_{1} & if $w = 0$,
    \cr 0_{2} & if $z = 0$. \quad $\diamondsuit$}$$
\end{example}

\begin{example}\label{nobasechange}
  {\it An AV--quotient without base--change property. }  Let $\Delta$
  be the fan in $\RR^{3}$ that has the maximal cones
  $$
  \sigma_{1} := \cone(-e_{1},e_{2}, e_{1}+e_{2}+e_{3}), \quad
  \sigma_{2} := \cone(e_{1},e_{2}, e_{1}+e_{2}+e_{3}). $$
  Consider the
  projection $P_{0} : \ZZ^{3} \to \ZZ^{2}$, $(u,v,w) \mapsto (u,v)$.
  Let $H$ be the subtorus of the acting torus of $X_{\Delta}$
  corresponding to the kernel of $P_{0}$. Then the TV--quotient $p
  \colon X_{\Delta} \to X_{\Delta} \tq H$ for the action of $H$ arises
  from the map $P \colon \ZZ^{3} \mapsto \ZZ$, $(u,v,w) \mapsto v$
  from $\Delta$ to the fan of faces of $\RR_{\ge 0}$.

  In particular, $X_{\Delta} \tq H = \CC$ is of dimension one.
  According to \ref{codimtwo}, $p$ is also an AV--quotient. But for
  the acting torus $\CC^{*}$ of $X_{\Delta} \tq H$, the open set
  $p^{-1}(\CC^{*})$ has, again by \ref{codimtwo}, a two--dimensional
  AV--quotient. \quad $\diamondsuit$
\end{example}

\section{A $\CC^{*}$--Action without AV--Quotient but with PV--Quotient}

We consider the open toric subvariety $X := \CC^{2} \times
(\CC^{*})^{2} \cup (\CC^{*})^{2} \times \CC^{2}$ of $\CC^{4}$ and the
action of the one--dimensional subtorus
$$H := \{(t,t,1,t^{-1}); \; t \in \CC^{*}\} \subset (\CC^{*})^{4}.$$

\begin{proposition}
\begin{enumerate}
\item There is a PV--quotient for the action of $H$ on $X$.
\item The action of $H$ on $X$ admits no AV--quotient.
\end{enumerate}
\end{proposition}

\proof Note that $X$ arises from the fan $\Delta$ in $\ZZ^{4}$ that has
$\sigma_{1} := \cone(e_{1},e_{2})$ and $\sigma_{2} :=
\cone(e_{3},e_{4})$ as its maximal cones. Let $F \colon \ZZ^{4} \to
\ZZ^{3}$ denote the lattice homomorphism defined by
$$
F(e_1) := e_1, \quad F(e_2) := e_2, \quad F(e_3) := e_3, \quad
F(e_4) := e_1+e_2.$$
To prove i), set $\tau_{i} := F_{\RR}(\sigma_{i})$ and define a
system $\S$ of fans in $\ZZ^{3}$ by $\Delta_{ii} :=
\mathfrak{F}(\tau_{i})$, where $i = 1,2$ and $\Delta_{12} :=
\Delta_{21} := \{ \{0 \}\}$.

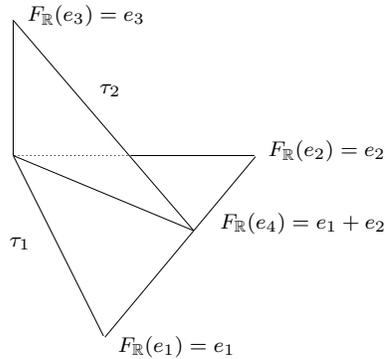
\begin{center}
  \begin{picture}(0,0)%
\includegraphics{bild1.pstex}%
\end{picture}%
\setlength{\unitlength}{1657sp}%
\begingroup\makeatletter\ifx\SetFigFont\undefined%
\gdef\SetFigFont#1#2#3#4#5{%
  \reset@font\fontsize{#1}{#2pt}%
  \fontfamily{#3}\fontseries{#4}\fontshape{#5}%
  \selectfont}%
\fi\endgroup%
\begin{picture}(3870,5190)(1756,-5461)
\put(3376,-5461){\makebox(0,0)[lb]{\smash{\SetFigFont{8}{9.6}{\rmdefault}{\mddefault}{\updefault}$F_\RR(e_1) = e_1$}}}
\put(5626,-2536){\makebox(0,0)[lb]{\smash{\SetFigFont{8}{9.6}{\rmdefault}{\mddefault}{\updefault}$F_\RR(e_2) = e_2$}}}
\put(4906,-3616){\makebox(0,0)[lb]{\smash{\SetFigFont{8}{9.6}{\rmdefault}{\mddefault}{\updefault}$F_\RR(e_4) = e_1+e_2$}}}
\put(1756,-3886){\makebox(0,0)[lb]{\smash{\SetFigFont{8}{9.6}{\rmdefault}{\mddefault}{\updefault}$\tau_1$}}}
\put(3106,-1591){\makebox(0,0)[lb]{\smash{\SetFigFont{8}{9.6}{\rmdefault}{\mddefault}{\updefault}$\tau_2$}}}
\put(2026,-511){\makebox(0,0)[lb]{\smash{\SetFigFont{8}{9.6}{\rmdefault}{\mddefault}{\updefault}$F_\RR(e_3) = e_3$}}}
\end{picture}

\end{center}

According to \cite{acha2}, Theorem 6.7, the toric morphism $X \to
X_{\S}$ defined by $F$ is a good prequotient for the action of $H$ on
$X$. In particular, it is a categorical prequotient.

We prove ii). Assume that there exists an AV--quotient $p \colon X \to Y$ for
the action of $H$ on $X$. We lead this to a contradiction by
presenting an $H$--invariant map that does not factor through
$p$. Consider
$$ f \colon X \to \CC^{3},  \qquad (x_{1},x_{2},x_{3},x_{4}) \mapsto
(x_{1}x_{4},x_{2}x_{4},x_{3}). $$
Note that
$$ f(X) = \CC^{3} \setminus (\{0\} \times \CC^* \times \{0\} \cup
\CC^* \times \{0\} \times \{0\}). $$
In particular, $f(X)$ is not open in $\CC^{3}$. We describe $f$ in
terms of fans. Let $\tau := \cone(e_{1},e_{2},e_{3}) \subset
\RR^{3}$. Then $f$ is just the toric morphism $X \to X_{\tau} =
\CC^{3}$ defined by the lattice homomorphism $F$. Thus it follows from
\cite{acha} that $f$ is the TV--quotient for the action of $H$ on $X$.

By its universal property, $p$ is surjective and there is a regular
map $\t{f} \colon Y \to X_{\tau}$ such that $f = \t{f} \circ p$. We
claim that all fibres of $\t{f}$ are of dimension zero. To see this let
$\varrho := \cone(e_{3}) \in \RR^{3}$ and note that surjectivity of
$p$ implies
$$ Y = \t{f}^{-1}(X_{\tau_{1}}) \cup \t{f}^{-1}(X_{\varrho}) \cup
\t{f}^{-1}(x_{\tau}). $$
By \cite{acha}, Example 3.1, the map $f_{1} :=
f\vert_{X_{\sigma_{1}}} \colon X_{\sigma_{1}} \to X_{\tau_{1}}$ is an
algebraic quotient for the action of $H$ on $X_{\sigma_{1}}$. Hence
we obtain a regular map $g \colon X_{\tau_{1}} \to Y$ and a
commutative diagram

$$\matrix{%
  X_{\sigma_1} & \subset & X & \topto{f} & X_{\tau} \cr
  \lmapdown{f_1} & & \rmapdown{p} & \ \rmapne{\t{f}} & \cr
X_{\tau_1} & \topto{g} & Y & &.  \cr}%
$$

Note that $\t{f} \circ g$ is necessarily an isomorphism and hence $g$
is an open embedding. Let $T := (\CC^{*})^{4}$ be the acting torus of
$p$ and set $\varrho_{4} := \cone(e_{4}) \subset
\RR^{4}$. By surjectivity of $p$ and the Fibre Formula 3.5 of
\cite{acha2}, we have 
$$ \t{f}^{-1}(X_{\tau_{1}}) = p(f^{-1}(X_{\tau_{1}})) =
p(X_{\sigma_{1}} \cup T \mal x_{\varrho_{4}}) = p(X_{\sigma_{1}}) =
g(X_{\tau_{1}}).$$
Here the third equality is a consequence of Lemma
\ref{xsigmaigleichewerte} ii). So $\t{f}$ is injective on
$\t{f}^{-1}(X_{\tau_{1}})$. A similar argument shows that $\t{f}$ is
injective on $\t{f}^{-1}(X_{\varrho})$. To verify the claim, we still
have to consider the fibre $\t{f}^{-1}(x_{\tau})$. Again by  Lemma
\ref{xsigmaigleichewerte} ii) one has
$$\t{f}^{-1}(x_{\tau}) = p(T \mal x_{\sigma_{2}}) \subset
p(\b{T \mal x_{\varrho_{4}}}) \subset \b{p(T \mal x_{\varrho_{4}})} =
\b{p(T \mal x_{\sigma_{1}})} = \b{\t{f}^{-1}(T_{1} \mal x_{\tau_{1}})}
.$$
Here $T_{1}$ denotes the acting torus of $X_{\tau}$. Since the closure
of $\t{f}^{-1}(T_{1} \mal x_{\tau})$ is contained in
$\t{f}^{-1}(\b{T_{1} \mal x_{\tau}})$, the above inclusion yields
$$ \t{f}^{-1}(\b{T_{1} \mal x_{\tau_{1}}}) = \t{f}^{-1}(T_{1} \mal
x_{\tau_{1}}) \cup \t{f}^{-1}(x_{\tau}) = \b{\t{f}^{-1}(T_{1} \mal
  x_{\tau_{1}})},$$
i.e., $\t{f}^{-1}(x_{\tau})$ is contained in the closure of of
$\t{f}^{-1}(T_{1} \mal x_{\tau_{1}})$. Moreover, we know that
$\t{f}^{-1}(T \mal x_{\tau_{1}}) = g(T \mal x_{\tau_{1}})$ is locally
closed of dimension one. Thus $\t{f}^{-1}(x_{\tau})$ is of dimension
zero and our claim is proved.

To conclude the proof, observe that by Zariski's Main Theorem, $\t{f}$
is an open embedding. This contradicts the fact that $f(X)$ is not
open in $X_{\tau}$. \endproof

\medskip

In fact the arguments used in our proof are chosen to work also in the
category of analytic spaces (for the existence of $g$ use
\cite{Sn}). Thus we obtain:

\begin{proposition}
The action of $H$ on $X$ does not admit categorical
quotients in the categories of analytic and algebraic
spaces. \endproof
\end{proposition}

\section{A $\CC^{*}$--Action admitting neither an AV--Quotient nor a
  PV--Quotient} 

Let $X$ denote the smooth four--dimensional toric variety obtained by
glueing the two affine charts $X_1=\CC^4$ and $X_2=\CC^3\times \CC^*$
along the common subset $(\CC\times\CC^*)^2$, using the glueing map
$$(t_1,t_2,t_3,t_4) \mapsto (t_1t_{2}^{2},t_2^{-1},t_3,t_4).$$ 
Let $T := (\CC^{*})^{4}$ denote the acting torus of $X$. 
We consider the action of the one--dimensional subtorus $H \subset T$
on $X$, where
$$ H := \{(t^{-2},1,t,t); \; t \in \CC^{*}\}.$$ 

\begin{proposition}\label{exconstr}
  There is neither an AV--quotient nor a PV--quotient for the action
  of $H$ on $X$.
\end{proposition}

\proof We will consider the TV--quotient $f \colon X \to X'$ for the
action of $H$ on $X$, which will turn out to be non--surjective. The
assumption that $f$ factors through a surjective $H$--invariant
regular map onto a complex prevariety will then lead to a
contradiction.

As before, we first describe the situation in terms of fans.
Let $e_1, \ldots, e_4$ denote the canonical basis vectors of $\RR^{4}$
and let $\Delta$ be the fan in $\RR^4$ with the maximal cones
$$
\sigma_1 := {\rm cone}(e_1, e_2, e_3, e_4) \quad
\text{and}\quad\sigma_2 := {\rm cone}(e_1, 2e_1 - e_2, e_3). $$
Note
that $\sigma_1 \cap \sigma_2$ is the cone spanned by $e_1$ and $e_3$.
The toric variety $X_{\Delta}$ associated to $\Delta$ equals $X$. In
order to describe $f$, consider the fan $\Delta'$ in $\RR^3$ with the
maximal cones
$$ \tau_1 := \cone(e_{1}-e_{2}, e_{1} + e_{3}, e_{1} - e_{3}), \quad
\tau_2 := \cone(e_{1}+e_{2}, e_{1}+e_{3}, e_{1}-e_{3}).$$

Then $f \colon X \to X'$ arises from the map $F \colon
\ZZ^4 \to \ZZ^3$ of the fans $\Delta$ and $\Delta'$ that, with respect
to the canonical bases, is given by the matrix

$$
\left[ \matrix{ 1 & 1 & 1 & 1 \cr 0 &-1 & 0 & 0 \cr 0 & 0 & 1 &-1
    \cr} \right]. $$

Note that $F_{\RR}(\sigma_1) = \tau_1$, whereas $F_{\RR}(\sigma_2)
\neq \tau_2$. More precisely, there is exactly one cone in $\Delta'$
whose relative interior does not intersect $F_{\RR}(\sigma_1) \cup
F(\sigma_2)$, namely the face $\tau$ of $\tau_2$ spanned by
$(1,1,0)$ and $(1,0,-1)$.

Let $T'$ denote the acting torus of $X_{\Delta'}$. Then we obtain
$f(X_\Delta) = X_{\Delta'} \setminus T' \mal x_{\tau}$. In particular,
$f$ is not surjective and $f(X_\Delta)$ is not open in $X_{\Delta'}$.

Now, assume that there is an AV-- or a PV--quotient
 for the action of $H$ on $X_{\Delta}$. Then, in both cases,
we have a surjective regular $H$--invariant map $p \colon X_\Delta \to Y$
onto a complex prevariety $Y$ and a regular map $\t{f} \colon Y \to
X_{\Delta'}$ such that the diagram $f = \t{f} \circ p$.

Note that $\t{f}$ is compatible with the induced (set
 theoretical) action of $T$ on $Y$, i.e., if $\varphi \colon T \to T'$
 denotes the homomorphism of the acting tori associated to $f$, then
 we have $\t{f}(t \mal y) = \varphi(t) \mal \t{f}(y)$ for all $t\in T$
 and $y\in Y$.

 We claim that $\t{f}$ has finite fibres and is injective over an open
 set of $X_{\Delta'}$. First consider the open affine
 toric subvariety $X_{\sigma_1}$ of $X_\Delta$. By \cite{acha},
 Example 3.1, the toric morphism $f_1 \colon X_{\sigma_1} \to X_{\tau_1}$
 defined by $F$ is the algebraic quotient for the action of $H$ on
 $X_{\sigma_1}$. Thus, there is a regular map $g \colon
 X_{\tau_1} \to Y$ such that the diagram

$$\matrix{%
  X_{\sigma_1} & \subset & X_\Delta & \topto{f} & X_{\Delta'} \cr
  \lmapdown{f_1} & & \rmapdown{p} & \ \rmapne{\t{f}} & \cr
X_{\tau_1} & \topto{g} & Y & &  \cr}%
$$

is commutative (see e.g. \cite{acha2}, Proposition 6.4). It follows
that $\t{f} \circ g$ defines an automorphism of $X_{\tau_1}$. Since
$p$ is surjective we have
$$
\t{f}^{-1}(X_{\tau_1}) = p(f^{-1}(X_{\tau_1})) = p(X_{\sigma_1}) =
g(X_{\tau_1}). $$
Consequently, $\t{f}$ is injective on the set
$\t{f}^{-1}(X_{\tau_1})$. Now consider $\sigma' := \cone(e_3,
2e_1-e_2)\in\Delta$ and set $\tau' := F_\RR(\sigma')$. By looking at
the toric morphism $f_2 \colon X_{\sigma'} \to X_{\tau'}$ induced by $f$, we
obtain with similar arguments as above that $\t{f}$ is injective over
$X_{\tau'}$.

\begin{center}
 \begin{picture}(0,0)%
\includegraphics{bild2.pstex}%
\end{picture}%
\setlength{\unitlength}{1657sp}%
\begingroup\makeatletter\ifx\SetFigFont\undefined%
\gdef\SetFigFont#1#2#3#4#5{%
  \reset@font\fontsize{#1}{#2pt}%
  \fontfamily{#3}\fontseries{#4}\fontshape{#5}%
  \selectfont}%
\fi\endgroup%
\begin{picture}(4950,4740)(1576,-3661)
\put(4051,-3661){\makebox(0,0)[lb]{\smash{\SetFigFont{8}{9.6}{\rmdefault}{\mddefault}{\updefault}$F_\RR(e_4)$}}}
\put(3601,-961){\makebox(0,0)[lb]{\smash{\SetFigFont{8}{9.6}{\rmdefault}{\mddefault}{\updefault}$F_\RR(\sigma_1)$}}}
\put(4726,-736){\makebox(0,0)[lb]{\smash{\SetFigFont{8}{9.6}{\rmdefault}{\mddefault}{\updefault}$F_\RR(\sigma_2)$}}}
\put(4051,-1636){\makebox(0,0)[lb]{\smash{\SetFigFont{8}{9.6}{\rmdefault}{\mddefault}{\updefault}$F_\RR(e_1)$}}}
\put(5401,-1636){\makebox(0,0)[lb]{\smash{\SetFigFont{8}{9.6}{\rmdefault}{\mddefault}{\updefault}$F_\RR(\sigma)$}}}
\put(1576,-1411){\makebox(0,0)[lb]{\smash{\SetFigFont{8}{9.6}{\rmdefault}{\mddefault}{\updefault}$F_\RR(e_2)$}}}
\put(5176,-61){\makebox(0,0)[lb]{\smash{\SetFigFont{8}{9.6}{\rmdefault}{\mddefault}{\updefault}$F_\RR(\sigma')$}}}
\put(6526,-1411){\makebox(0,0)[lb]{\smash{\SetFigFont{8}{9.6}{\rmdefault}{\mddefault}{\updefault}$F_\RR(2e_1-e_2)$}}}
\put(4051,839){\makebox(0,0)[lb]{\smash{\SetFigFont{8}{9.6}{\rmdefault}{\mddefault}{\updefault}$F_\RR(e_3)$}}}
\end{picture}

\end{center}
\begin{center}
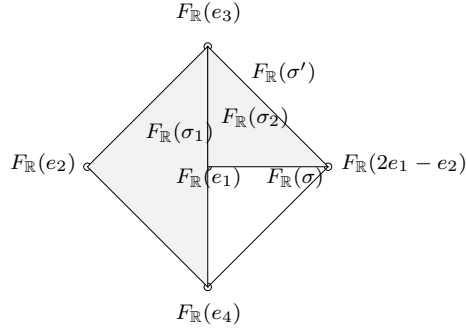
  {\small Intersection of $\Delta'$ with the plane
    defined by $x=1$ in $\RR^3$. }
\end{center}

Thus, to obtain our claim it remains to consider the fibre
$\t{f}^{-1}(x_{\tau_2})$. Note that according the Fibre Formula 3.5 of
\cite{acha2} one has

$$f^{-1}(x_{\tau_2}) = T \mal x_{\sigma_2} \cup T \mal x_{\sigma},$$
where $\sigma := \cone(e_1,2e_1 - e_2)\in\Delta$. Let $\varrho :=
\cone(e_{1}) \subset \RR^{4}$. We claim that $T \mal p(x_{\varrho})$
is locally closed of dimension one. This follows from the fact that
$T_{1} \mal f_{1}(x_{\varrho})$ and $T' \mal f(x_{\varrho})$ are locally
closed of dimension one.

Now note that $\varrho \prec \sigma_{2}$ and $\varrho \prec
\sigma$. Consequently $T \mal x_{\sigma_{2}}$ and $T \mal x_{\sigma}$
are contained in the closure of the orbit $T \mal x_{\varrho}$. This
implies
$$ T \mal p(x_{\sigma_{2}}) \cup T \mal p(x_{\sigma}) \subset \b{T
  \mal p(x_{\varrho})}. $$
Since $f(x_{\sigma_{2}}) \ne f(x_{\varrho}) \ne f(x_{\sigma})$, we
  obtain that the $T'$--orbits through $f(x_{\sigma_{2}})$ and
  $f(x_{\sigma})$ do not meet $T' \mal f(x_{\varrho})$. Thus we have
  even
$$ T \mal p(x_{\sigma_{2}}) \cup T \mal p(x_{\sigma}) \subset \b{T
  \mal p(x_{\varrho})} \setminus T \mal p(x_{\varrho}). $$
In other words, $\t{f}^{-1}(x_{\tau_{2}}) = p(f^{-1}(x_{\tau_{2}}))$
  consists of finitely many points. Thus we verified that $\t{f}$ has
  finite fibres.

Now, cover $Y$ by open affine charts $U_{1}, \ldots, U_{r}$ and set $U
:= T \mal p(x_{0}) = \t{f}^{-1}(T' \mal x'_{0})$. Then each
restriction $\t{f}_{i} := \t{f}\vert_{U_{i}}$ has finite fibres and is
injective along the non--empty open set $U_{i} \cap U$. Since $Y$ and
$X_{\Delta'}$ are normal, we obtain that the $\t{f}_{i}$ are open maps. This
yields openness of $f(X) = \bigcup \t{f}_{i}(U_{i})$ and we arrive at
a contradiction. \endproof

\section{A $\CC^*$--Action with AV--Quotient but without PV--Quotient}

We consider the smooth three--dimensional toric variety $X$ obtained
from glueing two copies of $\CC^{3}$ along the open subset $\CC\times
(\CC^*)^2$ by the following map:
$$(x_{1},x_{2},x_{3})\mapsto
(x_{1}x_{2}^{2}x_{3}^{2},x_{2}^{-1},x_{3}^{-1}).$$ 
In terms of convex geometry, $X$ is the toric variety arising
from the fan $\Delta$ in $\ZZ^{3}$ that has the maximal cones
$$\sigma_1  :=  \cone(e_{1},e_{1}-e_{2},e_{1}+e_{2}+e_{3}), \quad  
\sigma_2  :=  \cone(e_{1},e_{1}+e_{2},e_{1}-e_{2}-e_{3}).$$
Moreover, let $\t{X}$ denote the affine toric variety defined by the
fan $\t{\Delta}$ of faces of $\sigma := \cone(e_1+e_2,e_1-e_2)$ in
$\ZZ^2$.  Let $P \colon \ZZ^3 \to \ZZ^2, (x,y,z)\mapsto (x,y)$ denote
the projection.

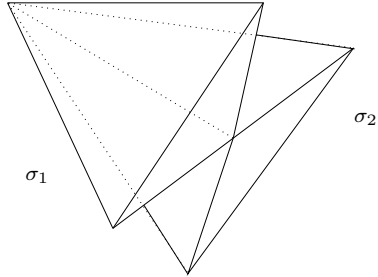
\begin{center}
  \begin{picture}(0,0)%
\includegraphics{bild3.pstex}%
\end{picture}%
\setlength{\unitlength}{1657sp}%
\begingroup\makeatletter\ifx\SetFigFont\undefined%
\gdef\SetFigFont#1#2#3#4#5{%
  \reset@font\fontsize{#1}{#2pt}%
  \fontfamily{#3}\fontseries{#4}\fontshape{#5}%
  \selectfont}%
\fi\endgroup%
\begin{picture}(5234,4109)(879,-3923)
\put(1171,-2491){\makebox(0,0)[lb]{\smash{\SetFigFont{8}{9.6}{\rmdefault}{\mddefault}{\updefault}$\sigma_1$}}}
\put(6091,-1621){\makebox(0,0)[lb]{\smash{\SetFigFont{8}{9.6}{\rmdefault}{\mddefault}{\updefault}$\sigma_2$}}}
\end{picture}

\end{center}

Then $P_{\RR}(\sigma_i) = \sigma$, so $P$ is a map of the fans $\Delta$ and
$\t{\Delta}$. Set $L := \ker(P)$. Note that $P$ is universal with
respect to $L$--invariant maps of fans (see \cite{acha}, Section 2) and
also with respect to $L$--invariant maps of systems of fans (see
\cite{acha2}, Section 7).

Now let $H$ be the one-dimensional subtorus of the acting torus $T =
(\CC^*)^3$ of $X$ that corresponds to $L$. Moreover, let $X_i \subset
X$ be the affine open subset corresponding to $\sigma_i$. Then
Proposition \ref{sepcrit}, \cite{acha2} Section 7 and \cite{acha},
Example 3.1 yield the following

\begin{proposition}
  The toric morphism $p \colon X \to \t{X}$ associated to $P$
  satisfies
\begin{enumerate}\label{properties}
\item $p$ is the AV--quotient for the action of $H$.
\item $p$ is the TP--quotient for the action of $H$.
\item The restriction $p_i := p \vert_{X_i} \colon X_i \to \t{X}$ is
  the algebraic quotient. \endproof
\end{enumerate}
\end{proposition}

But as we will show below, $p$ does not satisfy the universal property
of a PV--quotient in the category of arbitrary
prevarieties.  In fact, we even obtain

\begin{proposition}
  The action of $H$ on $X$ admits no PV--quotient.
\end{proposition}

\proof Assume that there is a PV--quotient $q \colon X \to
Y$. We claim that $Y$ is a toric prevariety and $q$ a toric morphism.
Note that there is an induced (set theoretical) $T$--action on $Y$
such that $q$ is equivariant. By the universal property of $q$ and
Proposition~\ref{properties}, there are commutative diagrams
$$\matrix{ X_i & \subset & X & \cr \lmapdown{p_i} & & \lmapdown{q} &
  \rmapse{p} \cr \t{X} & \botmapright{r_i} & Y & \botmapright{r} &
  \t{X} }
$$

of $T$--equivariant regular maps. Since $r \circ r_i = \id_{\t{X}}$
holds, each $r_i$ is injective, so Zariski's Main Theorem implies that
the $r_i$ are open embeddings. Since $X$ is covered by the $X_i$ and
$q$ is surjective, we obtain that $Y$ is covered by the $T$--stable
affine open subspaces
$$Y_i := r_i(\t{X}) = q(X_i).$$
In particular it follows that the
induced $T$--action on $Y$ is regular. So $Y$ is a toric prevariety
and $q$ is a toric morphism. This readily implies that $q$ satisfies
the universal property of a TP--quotient for the action of $H$ on
$X$. According to Remark \ref{properties}, we may assume $q = p$.

In order to show that $p$ is not a PV--quotient we
construct a map $f \colon X \to Z$ of prevarieties that does not
factor through $p$. Consider the maps $p_i$ defined above and the
distinguished points
$$
x_1 := x_{\varrho_1} \in X_1, \qquad x_2 := x_{\varrho_{2}} \in
X_2. $$
where $\varrho_{1} := \RR_{\ge 0} (e_1-e_2)$ and
$\varrho_{2} := \RR_{\ge 0} (e_1-e_2-e_{3})$. Note that the point $z
:= p(x_1) = p(x_2)$ does not lie in $p(X_1 \cap X_2)$. Consequently
the maps $p_i$ glue together to a regular map
$$
f \colon X = X_1 \cup_{X_1 \cap X_2} X_2 \to \t{X} \cup_{\t{X}
  \setminus \{z\}} \t{X} =: Z $$
of prevarieties. Since $f$ separates
the points $x_1$ and $x_2$, there is no set--theoretical factorization
of $f$ through $p$. \endproof

\bibliography{}

\bigskip

\bigskip

\noindent Fachbereich Mathematik und Statistik, \\
\noindent Universit\"at Konstanz \\
\noindent D--78457 Konstanz\\
\noindent E--mail: {\tt Annette.ACampo@uni-konstanz.de}\\
\noindent \hphantom{\hbox{E--mail:}} {\tt Juergen.Hausen@uni-konstanz.de}

\end{document}